\documentclass[12pt]{article}
\usepackage{amssymb, url}
\usepackage{amsmath,amsthm}
\usepackage{mathtools}
\usepackage{dsfont}
\usepackage{verbatim}
\usepackage{graphicx}
\usepackage{epstopdf}
\usepackage{fullpage}
\usepackage[hang,flushmargin]{footmisc}
\usepackage{enumitem}
\usepackage{tikz}

\usepackage{datetime}
\usepackage{bbm}

\usepackage{color}

\newtheorem{theorem}{Theorem}[]

\newtheorem{lemma}{Lemma}
\newtheorem{conjecture}{Conjecture}

\theoremstyle{definition}

\newtheorem{claim}{Claim}[theorem]

\providecommand{\keywords}[1]
{
  \small	
  \textbf{\textit{Keywords---}} #1
}

\begin{document}
\title{Maximally Edge-Connected Realizations and Kundu's $k$-factor Theorem}
\author{James M.\ Shook$^{3,4}$}

\footnotetext[3]{National Institute of Standards and Technology, Computer Security Division, Gaithersburg, MD; {\tt james.shook@nist.gov}.}
\footnotetext[4]{Official Contribution of the National Institute of Standards and Technology; Not subject to copyright in the United States.}
% \today

\maketitle
\begin{abstract}
    A simple graph $G$ with edge-connectivity $\lambda(G)$ and minimum degree $\delta(G)$ is maximally edge connected if $\lambda(G)=\delta(G)$. In 1964, given a non-increasing degree sequence $\pi=(d_{1},\ldots,d_{n})$, Jack Edmonds showed that there is a realization $G$ of $\pi$ that is $k$-edge-connected if and only if $d_{n}\geq k$ with  $\sum_{i=1}^{n}d_{i}\geq 2(n-1)$ when $d_{n}=1$. We strengthen Edmonds's result by showing that given a realization $G_{0}$ of $\pi$ if $Z_{0}$ is a spanning subgraph of $G_{0}$ with $\delta(Z_{0})\geq 1$ such that $|E(Z_{0})|\geq n-1$ when $\delta(G_{0})=1$, then there is a maximally edge-connected realization of $\pi$ with $G_{0}-E(Z_{0})$ as a subgraph. Our theorem tells us that there is a maximally edge-connected realization of $\pi$ that differs from $G_{0}$ by at most $n-1$ edges. For $\delta(G_{0})\geq 2$, if $G_{0}$ has a spanning forest with $c$ components, then our theorem says there is a maximally edge-connected realization that differs from $G_{0}$ by at most $n-c$ edges. As an application we combine our work with Kundu's $k$-factor Theorem to show there is a maximally edge-connected realization with a $(k_{1},\dots,k_{n})$-factor for $k\leq k_{i}\leq k+1$ and present a partial result to a conjecture that strengthens the regular case of Kundu's $k$-factor theorem.
\end{abstract}


\section{Introduction}
We only consider simple graphs, and see Diestel \cite{Diestel2016} for terminology not defined here. For a graph $G=(V,E)$ and $v\in V$, we let $deg_{G}(v)$ denote the number of neighbors of $v$ in $G$, and let $\Delta(G)$ and $\delta(G)$ denote the maximal and minimal degrees of $G$, respectively. We let $\lambda(G)$ denote edge-connectivity, and say $G$ is maximally edge-connected if $\lambda(G)=\delta(G)$.  

A sequence of non-negative integers $\pi=(d_{1},\ldots,d_{n})$ is a degree sequence if there exists a graph $G$ with vertex set $V=\{v_{1},\ldots,v_{n}\}$ such that $deg_{G}(v_{i})=d_{i}$. Such a graph $G$ is said to realize or is a realization of $\pi$. We call the sequence $\pi$ graphic if it is a degree sequence. If instead a graph $G$ with vertex set $V=\{v_{1},\ldots,v_{n}\}$ is given, then we let $\pi(G)=(deg_{G}(v_{1}),\ldots, deg_{G}(v_{n}))$. If every entry of a degree sequence is non-zero, then we say the degree sequence is positive. We let $\mathcal{R}(\pi)$ denote the set of realizations of the graphic sequence $\pi$. For a graph $F$, with vertex set $V$, we let $\mathcal{R}(\pi, F)\subseteq \mathcal{R}(\pi)$ be the set of all realizations whose set of edges include $E(F)$. We have $\mathcal{R}(\pi)=\mathcal{R}(\pi,\emptyset)$, and we write $\mathcal{R}(G,F)$ for $\mathcal{R}(\pi(G),F)$.

The conditions for when a graph is maximally edge-connected are well studied \cite{Hellwig2008}. In particular, given a graph $G$,  Bollab\'{a}s  \cite{Bollobas1979} and then extended upon in \cite{Dankelmann2000,Hellwig2003} found degree conditions for $\pi(G)$ that imply $G$ is maximally edge-connected. On the other hand, Jack Edmonds \cite{Edmonds1964} gave necessary and sufficient conditions for when a degree sequence has a realization that is maximally edge-connected.

\begin{theorem}[\cite{Edmonds1964}]\label{thm:Edmonds}
For a non-increasing degree sequence $\pi=(d_{1},\ldots,d_{n})$, there is a $G\in \mathcal{R}(\pi)$ that is $k$-edge-connected if and only if $d_{n}\geq k$ with  $\sum_{i=1}^{n}d_{i}\geq 2(n-1)$ when $d_{n}=1$. 
\end{theorem}
Edmonds did this work in 1964 while at the National Institute of Standards and Technology (NIST was known as the National Bureau of Standards (NBS) when Edmonds did his work.), and besides a short constructive proof of Theorem~\ref{thm:Edmonds} provided by Kleitman and Wang \cite{Wang1974} in 1974, only within the last ten years has Theorem~\ref{thm:Edmonds} been extended. Gu and Lai \cite{Gu2013} generalized Theorem~\ref{thm:Edmonds} to $k$-edge-connected uniform hypergraphs, and around the same time Tian, Meng, Lai, and Zhang \cite{Tian2014} gave necessary and sufficient conditions for when a degree sequence has a realization that is super edge-connected. A super edge-connected graph is one where every minimum edge cut isolates a vertex with minimum degree. While studying edge-disjoint perfect matchings we needed to strengthen Theorem~\ref{thm:Edmonds} so that the realization preserved a subgraph found in the original graph. However, we are able to prove more and so we present that here.

\begin{theorem}\label{thm:connectedLeftoverFull}If there is a graph $G_{0}=(V,E)$ with edge-disjoint spanning subgraphs $F$ and $Z_{0}$ with $\delta(Z_{0})>\Delta(F)$ such that $|E(Z_{0})|\geq |V|-1$ when $\delta(G_{0})=1$, then there is a $G\in \mathcal{R}(G_{0},G_{0}-E(Z_{0}))$ such that $G-E(F)$ is maximally edge-connected.
\end{theorem}

It is not hard to derive Theorem~\ref{thm:Edmonds} from Theorem~\ref{thm:connectedLeftoverFull}. The first part of Theorem~\ref{thm:Edmonds} is trivial and the second part follows from Theorem~\ref{thm:connectedLeftoverFull} by letting $F$ be empty and $Z_{0}=G_{0}$. 

When $F$ is empty Theorem~\ref{thm:connectedLeftoverFull} says any graph $G_{0}$ differs from a maximally edge-connected realization of $\pi(G_{0})$ by at most $|V(G_{0})|-1$ edges. For $\delta(G_{0})\geq 2$, if we let $Z_{0}$ be a spanning forest with $c$ components, then we can show $G_{0}$ differs from some maximally edge-connected realization of $\pi(G_{0})$ by at most $|V(G_{0})|-c$ edges. Interestingly, this means if $G_{0}$ has a perfect matching, then there is some maximally edge-connected realization that differs form $G_{0}$ by at most $|V(G_{0})|/2$ edges.

Observe in Theorem~\ref{thm:connectedLeftoverFull} that if $F$ is maximally edge-connected or $F=\emptyset$, then $G$ is maximally edge connected. Moreover, if $F=\emptyset$ and $Z_{0}=G_{0}-E(H)$ where $H$ is a subgraph of $G_{0}$ such that $\Delta(H)\leq d_{n}-1$ and $|E(Z_{0})|\geq n-1$, then Theorem~\ref{thm:connectedLeftoverFull} says $H$ is a subgraph of some maximally edge-connected realization of $\pi$. With a simpler proof than Theorem~\ref{thm:connectedLeftoverFull} we may allow $\delta(Z_{0})\geq \Delta(F)$ when $Z_{0}=G_{0}-E(F)$ at the expense of lowering the edge connectivity of $G-E(F)$ by one when $\delta(G-E(F))$ is odd. 

\begin{theorem}\label{thm:simplifyConnectedLeftoverFull}If there is a subgraph $F$ of a graph $G_{0}=(V,E)$ with $\delta(G_{0}-E(F))\geq \Delta(F)$ such that $|E(G_{0}-E(F))|\geq |V|-1$ when $\delta(G_{0})=1$, then there is a $G\in \mathcal{R}(G_{0},F)$ such that $G-E(F)$ is maximally edge-connected when $\delta(G-E(F))$ is even and $(\delta(G-E(F))-1)$-edge-connected when $\delta(G-E(F))$ is odd.
\end{theorem}

\subsection{Kundu's $k$-factor Theorem}
For an application of Theorem~\ref{thm:connectedLeftoverFull} and Theorem~\ref{thm:simplifyConnectedLeftoverFull},  we look to Kundu's $k$-factor Theorem for inspiration. Recall a graph $G$ is said to have a $(k_{1},\dots,k_{n})$-factor if $G$ has a spanning subgraph with degree sequence $(k_{1},\dots, k_{n})$. If $k_{i}=k$ for all $i$, then we simply call the spanning subgraph a $k$-factor.
\begin{theorem}[Kundu's $k$-factor Theorem \cite{Kundu1974}]\label{kundu}
For $k\geq 0$, if $\pi=(d_{1},\ldots,d_{n})$ and $(d_{1}-k_{1},\ldots,d_{n}-k_{n})$ are both graphic such that $k\leq k_{i}\leq k+1$ for $1\leq i\leq n$, then there exists a realization of $\pi$ that has a $(k_{1},\dots,k_{n})$-factor.
\end{theorem}

By requiring each term of the degree sequence to be at least two and $k\geq 1$ we can use Theorem~\ref{thm:connectedLeftoverFull} to strengthen Kundu's Theorem. 

\begin{theorem}For $k\geq 1$, if $\pi=(d_{1},\ldots,d_{n})$ and $(d_{1}-k_{1},\ldots,d_{1}-k_{n})$ are both graphic such that $d_{i}\geq 2$ and $k\leq k_{i}\leq k+1$ for all $i\in \{1,\ldots, n\}$, then there exists a maximally-edge connected realization of $\pi$ that has a $(k_{1},\dots,k_{n})$-factor.
\begin{proof}By Kundu's $k$-factor Theorem, there exist a $G\in \mathcal{R}(\pi)$ that has a $(k_{1},\dots,k_{n})$-factor $H$. If we let $F = \emptyset$ and
$Z_{0} = H$, then Theorem~\ref{thm:connectedLeftoverFull} says there exists a maximally edge-connected realization $G'\in R(\pi)$ that contains the subgraph $G-E(H)$. Moreover, the edges of $G'$ not in the subgraph $G-E(H)$ form a $(k_{1},\dots, k_{n})$-factor.
\end{proof}
\end{theorem}

For a sequence $\pi=(d_{1},\ldots,d_{n})$, we let $\mathcal{D}_{k}(\pi)$ denote the sequence $(d_{1}-k,\ldots,d_{n}-k)$. Busch, Ferrara, Hartke, Jacobson, Kaul, and West \cite{Busch2012} showed that if $n$ is even and both $\pi$ and $\mathcal{D}_{k}(\pi)$ are graphic, then for $r\leq \min\{2,k\}$ there is a realization of $\pi$ with a $k$-factor that has $r$ edge-disjoint $1$-factors. Seacrest \cite{Seacrest2021} improved this for $r\leq \min\{4,k\}$. The work of \cite{Busch2012} and \cite{Seacrest2021} on edge-disjoint $1$-factors suggests a further strengthening of Kundu's theorem. 
\begin{conjecture}[\cite{Brualdi1978} and later in \cite{Busch2012}]\label{conj:KunduExpansion}Some realization of a degree sequence $(d_{1},\ldots,d_{n})$ with even $n$ has $k$ edge-disjoint $1$-factors if and only if $(d_{1}-k,\ldots,d_{n}-k)$ is graphic.
\end{conjecture}

Conjecture~\ref{conj:KunduExpansion} was first posed by Brualdi \cite{Brualdi1978} and then independently by Busch et al.\ in \cite{Busch2012}. The work of Seacrest \cite{Seacrest2021} shows the conjecture is true for $k\leq 5$. Busch et al.\ showed the conjecture holds for $d_{n}\geq \frac{n}{2}+k-2$ and $d_{1}\leq \frac{n}{2}+1$. In this paper we focus on large $k$.

A $1$-factorization of a graph is the partition of its edges into $1$-factors. Chetwyn and Hilton \cite{Chetwynd1985} described a $1$-factorization conjecture that says every $k$-regular graph with $k\geq 2\lceil \frac{n}{4}\rceil-1$ has a $1$-factorization. Impressively, Csaba, K\"{u}hn, Lo, Osthos, and Treglown affirmed this conjecture for $n$ sufficiently large. 

\begin{theorem}[\cite{Csaba2016}]\label{thm:Csaba}
There exists an $n_{0}\in \mathbb{N}$ such that the following holds. Let
$n,k \in \mathbb{N}$ be such that $n\geq n_{0}$ is even and $k\geq 2\lceil \frac{n}{4}\rceil-1$. Then every $k$-regular graph $G$ on $n$ vertices can be decomposed into $k$ edge-disjoint $1$-factors.
\end{theorem}

Thus, Theorem~\ref{thm:Csaba} says Conjecture~\ref{conj:KunduExpansion} is true for large $k$ and $n$ sufficiently large. However, we can say more now that we have Theorem~\ref{thm:simplifyConnectedLeftoverFull}.

The following classic result of Berge \cite{Berge1958} was expanded upon in \cite{Bollobas1985, Katerinis1993, Plesnik1972, Shiu2008}.

\begin{theorem}[\cite{Berge1958}]\label{thm:berge} All even ordered $(k-1)$-edge-connected $k$-regular graphs have a $1$-factor.
\end{theorem}
We use the result of Berge to show that for large $k$ we can find a $k$-factor with many edge-disjoint $1$-factors.

\begin{theorem}\label{thm:largek}Let $\pi=(d_{1},\ldots,d_{n})$ be a non-increasing degree sequence with even $n$ such that $\mathcal{D}_{k}(\pi)$ is graphic. If $k\geq \frac{d_{1}}{2}+r$ or $k\geq n-1-d_{n}+2r$, then $\pi$ has a realization with a $k$-factor that has at least $r+1$ edge-disjoint $1$-factors.
\begin{proof}We will first prove the case $k\geq d_{1}/2+r$. By Kundu's $k$-factor theorem there is some realization of $\pi$ with a $k$-factor. Let $i\leq k$ be the largest non-negative integer such that there is a $G_{i}\in \mathcal{R}(\pi)$ with a $(k-i)$-factor $H_{i}$ such that $F_{i}=G_{i}-E(H_{i})$ has $i$ edge-disjoint $1$-factors. Since $\delta(H_{i})=k-i\geq d_{1}/2+r-i\geq (k-i+\Delta(F_{i}))/2+r-i$, we see that $\delta(H_{i})\geq \Delta(F_{i})+2(r-i)$. Therefore, if $r\leq i$, then we may apply Theorem~\ref{thm:simplifyConnectedLeftoverFull} to find a $G_{i+1}\in \mathcal{R}(\pi,F_{i})$ such that $H_{i+1}=G_{i+1}-E(F_{i})$ is a $(k-i-1)$-edge-connected $(k-i)$-factor. However, we deduce a contradiction since Theorem~\ref{thm:berge} implies $H_{i+1}$ has a $1$-factor, and therefore, $G_{i+1}$ has a $k$-factor with $i+1$ edge-disjoint $1$-factors. Thus, $i\geq r+1$.

The case $k\geq n-1-d_{n}+2r$ can be proved with an application of the first part of this theorem to the non-increasing degree sequences  $(n-1-d_{n}+k,\ldots, n-1-d_{1}+k)$ and $(n-1-d_{n},\ldots, n-1-d_{1})$, the reverse order of the complements of $\mathcal{D}_{k}(\pi)$ and $\pi$, to find a realization with a $k$-factor that has $r+1$ edge-disjoint $1$-factors. The $k$-factor can then be mapped to a realization of $\pi$
\end{proof}
\end{theorem}

It would be interesting to see if Theorem~\ref{thm:berge} or results like it can be used to find more edge-disjoint $1$-factors. Considering the $1$-factorization conjecture and the work of Csaba et al., it seems likely that as $k$ increases the number of edge-disjoint $1$-factors in a maximally edge-connected regular graph increases. Thomassen in \cite{Thomassen2020} showed that the edges of a $k$-regular $k$-edge-connected graphs with some restrictions can be partitioned in various ways. Although, some of those restrictions may not be necessary. Thomassen, in the same paper, posed some nice conjectures and problems that avoids them. However, Mattiolo \cite{Mattiolo2022}, answering Problem~1 in \cite{Thomassen2020}, presented $k$-regular $k$-edge-connected graphs that cannot be partitioned into a $2$-factor and $k-2$ $1$-factors. Thus, our strategy maybe limited to finding many edge-disjoint $1$-factors, but not $k$ of them.

\section{Proofs}
Let $G=(V,E)$ be a graph with $X,Y\subseteq V$. We denote $G[X]$ to be the induced graph on $X$, and let $\overline{X}=V-X$. We denote $E_{G}(X,Y)$ to be the set of all edges of $G$ that have one end in $X$ and the other in $Y$, and let $e_{G}(X,Y)=|E_{G}(X,Y)|$ and $e_{G}(x,Y)=e_{G}(\{x\},Y)$.  We denote $\Gamma_{G}(X)$ to be the set of all vertices in $X$ that are adjacent in $G$ to vertices in $\overline{X}$. 

Let $G=(V,E)$ be a graph with $\lambda(G)<\delta(G)$, and let $A\subset V$. If $e_{G}(A,\overline{A})<\delta(G)$, then we say $A$ is weak. If $e_{G}(A,\overline{A})=\lambda(G)$, then we say $A$ is minimally weak. If $A$ is weak and $e_{G}(S,\overline{S})\geq \delta(G)$ for every $S\subset A$, then we say $A$ is critically weak.

The next two lemmas play an important role in the proofs of Theorem~\ref{thm:connectedLeftoverFull} and Theorem~\ref{thm:simplifyConnectedLeftoverFull}.

\begin{lemma}\label{lem:criticalSet}For a graph $G=(V,E)$ with $\lambda(G)<\delta(G)$, if $S\subseteq A\subset V$ such that $A$ is critically weak and $e_{G}(S,\overline{S})\leq \delta(G)$, then for any $X\subset V$ such that $X\cap S\neq \emptyset$ and $\overline{X}\cap S\neq \emptyset$ we see that  \begin{equation}\label{eq:criticalSet0}e_{G}(X\cap S, \overline{X}\cap S)\geq \Bigg\lceil\frac{\delta(G)}{2}\Bigg\rceil\end{equation} when $S\subset A$ and for $A=S$, we see that \begin{equation}\label{eq:criticalSet}e_{G}(X\cap A, \overline{X}\cap A)\geq \Bigg\lceil\frac{\delta(G)+1}{2}\Bigg\rceil
\end{equation} where equality in (\ref{eq:criticalSet}) implies \begin{equation}\label{eq:criticalSet3}\min\{e_{G}(X\cap A, \overline{A}),e_{G}(\overline{X}\cap A, \overline{A})\}\geq \Bigg\lfloor\frac{\delta(G)-1}{2}\Bigg\rfloor.
\end{equation}

\begin{proof}By the definition of a critically weak set we know that $e_{G}(S\cap X,\overline{S\cap X})\geq \delta(G)$ and $e_{G}(S\cap \overline{X},\overline{S\cap \overline{X}})\geq \delta(G)$. Thus, 
\begin{align}
    e_{G}(S,\overline{S})&=e_{G}(X\cap S,\overline{X\cap S})+e_{G}(\overline{X}\cap S,\overline{\overline{X}\cap S})-2e_{G}(X\cap S, \overline{X}\cap S)\notag\\
    &\geq 2(\delta(G)-e_{G}(X\cap S, \overline{X}\cap S)).\label{eq:1}
\end{align}
When $S\subset A$ we have by our condition on $S$ that $e_{G}(S,\overline{S})=\delta(G)$. Thus, solving (\ref{eq:1}) for $e_{G}(X\cap S, \overline{X}\cap S)$ we establish (\ref{eq:criticalSet0}). Using the fact that $e_{G}(A,\overline{A})\leq \delta(G)-1$ when $S=A$ we can establish (\ref{eq:criticalSet}) by solving (\ref{eq:1}) for $e_{G}(X\cap A, \overline{X}\cap A)$. Since $A$ is critically weak, when equality holds in (\ref{eq:criticalSet}) we see that 
\begin{align}
    \min\{e_{G}(X\cap A, \overline{A}),e_{G}(\overline{X}\cap A, \overline{A})\}&\geq \delta(G)-e_{G}(X\cap A, \overline{X}\cap A)\notag\\
    &= \delta(G)-\Bigg\lceil\frac{\delta(G)+1}{2}\Bigg\rceil=\Bigg\lfloor\frac{\delta(G)-1}{2}\Bigg\rfloor.\notag \qedhere
\end{align}
\end{proof}
\end{lemma}

The following lemma was proved in different forms in \cite{Dankelmann2000, Goldsmith1978}. We present our proof of the form we need.

\begin{lemma}\label{lem:weakBoundEdge}
For a graph $G=(V,E)$ with $\lambda(G)<\delta(G)$, if $A\subset V$ is weak, then $A-\Gamma_{G}(A)\neq \emptyset$ and therefore, every $x\in A-\Gamma_{G}(A)$ has a neighbor in $A-\Gamma_{G}(A)$.
\begin{proof}Suppose $\Gamma_{G}(A)=A$. In this case every vertex in $A$ must be adjacent in $G$ to at least $\delta(G)-(|A|-1)$ vertices in $\overline{A}$. Thus, \[\delta(G)-1\geq e_{G}(A,\overline{A})\geq |A|(\delta(G)-(|A|-1)).\] If we combine all terms of the last inequality onto the right hand side and simplify, then we see the contradiction \[0\geq (\delta(G)-|A|)(|A|-1)+1\] since $\delta(G)-1\geq |A|\geq 1$. Thus, $A-\Gamma_{G}(A)$ is not empty, and consequently, $N_{G}(x)\subseteq A$ for every $x\in A-\Gamma_{G}(A)$. Since $|\Gamma_{G}(A)|\leq \delta(G)-1$, every $x\in A-\Gamma_{G}(A)$ has a neighbor in $A-\Gamma_{G}(A)$.
\end{proof}
\end{lemma}

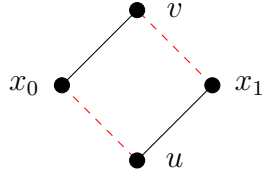
\begin{figure}
\centering
\begin{tikzpicture}
%edges
\draw (1,3) -- (2,4);
\draw (3,3) -- (2,2);
\draw[dashed, color=red] (1,3) -- (2,2);
\draw[dashed, color=red] (3,3) -- (2,4);
%nodes
\draw[fill=black] (1,3) circle (.1cm);
\draw[fill=black] (3,3) circle (.1cm);
\draw[fill=black] (2,2) circle (.1cm);
\draw[fill=black] (2,4) circle (.1cm);
%labels
\node[] at (.5,3) {$x_{0}$};
\node[] at (3.5,3) {$x_{1}$};
\node[] at (2.5,4) {$v$};
\node[] at (2.5,2) {$u$};
\end{tikzpicture}
\caption{An edge-exchange between edges $vx_{0}$ and $x_{1}u$ of $G$ and non-edges $x_{0}u$ and $vx_{1}$.}
\label{fig:Exchanges}
\end{figure}

We will often make use of an edge-exchange. This operation, see Figure~\ref{fig:Exchanges}, consists of exchanging two edges $vx_{0}$ and $x_{1}u$ from a graph $G$ with two edges $x_{0}u$ and $vx_{1}$ from $\overline{G}$ to create another realization of $\pi(G)$ while leaving all other edges the same. We will often refer to edges of $\overline{G}$ as non-edges of $G$. 

\subsection{Proof of Theorem~\ref{thm:simplifyConnectedLeftoverFull}}

Since the proof of Theorem~\ref{thm:simplifyConnectedLeftoverFull} follows a simplified proof of Theorem~\ref{thm:connectedLeftoverFull}, we present it here first.

\begin{proof}
We choose a $G\in \mathcal{R}(G_{0},F)$ such that 
\begin{enumerate}[label=(C\arabic*),ref=(C\arabic*)]
    \item\label{SimpconnectedLeftoverFull_choice:2} $\lambda(G-E(F))$ is maximized, and
    \item \label{SimpconnectedLeftoverFull_choice:3} subject to \ref{SimpconnectedLeftoverFull_choice:2}, we minimize the number of minimally weak sets in $G-E(F)$.
\end{enumerate} 
Let $Z=G-E(F)$. For a contradiction, we assume $\lambda(Z)<\delta(Z)$ and $\lambda(Z)<\delta(Z)-1$ when $\delta(Z)$ is odd.  As a consequence, we may choose an arbitrary minimally weak set $A_{0}\subseteq V(G)$ and critically weak sets $A\subseteq A_{0}$ and $B\subseteq \overline{A_{0}}$.

We choose an $a\in A-\Gamma_{Z}(A)$ and a $b\in B-\Gamma_{Z}(b)$ such that we give priority to an adjacent pair in $F$. If $a$ and $b$ are adjacent, then $N_{Z}(a)-N_{F}(b)\neq \emptyset$ and $N_{Z}(b)-N_{F}(a)\neq \emptyset$ since $\delta(Z)\geq \Delta(F)$. If $a$ and $b$ are not adjacent, then by our choice of $a$ and $b$, we deduce that $N_{F}(a)\cap (B-\Gamma_{Z}(b))=\emptyset$ and $N_{F}(b)\cap (B-\Gamma_{Z}(a))=\emptyset$. In either case we may choose an $a'\in N_{Z}(a)-N_{F}(b)$ and a $b'\in N_{Z}(b)-N_{F}(a)$.

Let $W$ be the realization of $\mathcal{R}(G_{0},F)$ created by exchanging the edges $aa'$ and $bb'$ of $Z$ with the non-edges $ab'$ and $ba'$ of $G$. By \ref{SimpconnectedLeftoverFull_choice:2} we know that $\lambda(Z)\geq \lambda(W-E(F))$. Let us first examine $A_{0}$ in $W-E(F)$. Since the edges $aa'$ and $bb'$ are not in $E_{Z}(A_{0},\overline{A_{0}})$ and the edges $ab'$ and $ba'$ are in $E_{W-E(F)}(A_{0},\overline{A_{0}})$, we see that $e_{W-E(F)}(A_{0},\overline{A_{0}})=e_{Z}(A_{0},\overline{A_{0}})+2$. Since $A_{0}$ is minimally weak in $Z$, we see that \[e_{W-E(F)}(A_{0},\overline{A_{0}})=e_{Z}(A_{0},\overline{A_{0}})+2=\lambda(Z)+2\geq \lambda(W-E(F))+2 .\] Therefore, $A_{0}$ is not minimally weak in $W-E(F)$.

We choose an arbitrary $X\subseteq V$ that is minimally weak in $W-E(F)$. Thus, \[\lambda(Z)\geq \lambda(W-E(F))=e_{W-E(F)}(X,\overline{X}).\] Suppose at most one of $aa'$ or $bb'$ is in $E_{Z}(X,\overline{X})$. We have three cases to consider. If both $aa'$ and $bb'$ are in $E_{Z}(X,X)$ or in $E_{Z}(\overline{X},\overline{X})$, then both $ab'$ and $ba'$ are in $E_{W-E(F)}(X,X)$ or in $E_{W-E(F)}(\overline{X},\overline{X})$. If one of $aa'$ or $bb'$ is in $E_{Z}(X,\overline{X})$, then one of $ab'$ or $ba'$ is in $E_{W-E(F)}(X,\overline{X})$. Thus, $e_{W-E(F)}(X,\overline{X})= e_{Z}(X,\overline{X})$ in the first two cases. Finally, if one of $aa'$ or $bb'$ is in $E_{Z}(X,X)$ and the other is in $E_{Z}(\overline{X},\overline{X})$, then both $ab'$ and $ba'$ are in $E_{W-E(F)}(X,\overline{X})$. Thus, $e_{W-E(F)}(X,\overline{X})\geq e_{Z}(X,\overline{X})$ in all three cases. From this we may deduce that \[\lambda(Z)\geq \lambda(W-E(F))=e_{W-E(F)}(X,\overline{X})\geq e_{Z}(X,\overline{X})\geq \lambda(Z).\] This implies  $\lambda(Z)=\lambda(W-E(X))$ and $X$ is minimally weak in $Z$. Thus, $W-E(F)$ satisfies \ref{SimpconnectedLeftoverFull_choice:2}, and since $X$ is an arbitrarily chosen minimally weak set of $W-E(F)$, we may conclude that every minimally weak set in $W-E(F)$ is minimally weak in $Z$.  However, since $A_{0}$ is not minimally weak in $W-E(F)$, we may, in contradiction with \ref{SimpconnectedLeftoverFull_choice:3}, conclude that $W-E(F)$ has fewer minimally weak sets than $Z$. Thus, $W$ contradicts our choice of $Z$, and therefore, we may assume $\{aa',bb'\}\subseteq E_{Z}(X,\overline{X})$. Since the edge-exchange between $aa'$ and $bb'$ only affects the two edges, we may conclude that $e_{W-E(F)}(X,\overline{X})\geq e_{Z}(X,\overline{X})-2$. By Lemma~\ref{lem:criticalSet} $e_{Z}(X\cap A, \overline{X}\cap A)$ and $e_{Z}(X\cap B, \overline{X}\cap B)$ are at least $\Big\lceil\frac{\delta(Z)+1}{2}\Big\rceil$. Thus, \[\lambda(Z)\geq e_{W-E(F)}(X,\overline{X})\geq e_{Z}(X,\overline{X})-2\geq 2\Bigg\lceil\frac{\delta(Z)+1}{2}\Bigg\rceil-2.\] However, this presents a contradiction since the right hand side of the last inequality is at least $\delta(Z)$ when $\delta(Z)$ is even and at least $\delta(Z)-1$ when $\delta(Z)$ is odd. Thus, $G$ satisfies Theorem~\ref{thm:simplifyConnectedLeftoverFull}.\qedhere
\end{proof}

\subsection{Proof of Theorem~\ref{thm:connectedLeftoverFull}}
Edmonds established Theorem~\ref{thm:Edmonds} by directly proving the $\delta(G-E(F))=1$ case and then used a strategy of reducing the number of weak sets when $\delta(G-E(F))\geq 2$. We follow the same strategy, but our job is more difficult since we have fewer edges with which we may exchange. We tackle this difficulty by carefully selecting critically weak sets and vertices so we may find edges with useful properties.

\begin{proof}
We choose a $G\in \mathcal{R}(G_{0},G_{0}-E(Z_{0}))$ such that 
\begin{enumerate}[label=(C\arabic*),ref=(C\arabic*)]
    \item\label{connectedLeftoverFull_choice:2} $\lambda(G-E(F))$ is maximized, and
    \item \label{connectedLeftoverFull_choice:3} subject to \ref{connectedLeftoverFull_choice:2}, we minimize the number of minimally weak sets in $G-E(F)$.
\end{enumerate} 
Let $H=G-E(F)$, and by contradiction we assume $\lambda(H)<\delta(H)$.  We let $Z=G-E(G_{0}-E(Z_{0}))$, and observe that $Z\in \mathcal{R}(Z_{0})$. 

\begin{claim}\label{cl:basicCase}$H$ is connected, and $\delta(H)\geq 2$.
\begin{proof}For a contradiction, we assume $H$ can be partitioned into the components $C_{1},\ldots,C_{t}$. Suppose there is a component $C_{l}$ that has a cycle containing an edge $aa'\in E(Z)$. Since $\delta(Z)>\Delta(F)$, we can choose an edge $bb'$ of $Z$ in some other component $C_{j}$ such that $b\notin N_{F}(a')$ and $a\notin N_{F}(b')$. We exchange the edges $aa'$ and $bb'$ of $Z$ with the non-edges $ab'$ and $ba'$ of $G$ to create new realizations $G'\in \mathcal{R}(G_{0},G_{0}-E(Z_{0}))$. Since $aa'$ was in a cycle of $C_{l}$ the vertices in $V(C_{l})\cup V(C_{i})$ form a component of $G'$, and therefore, $G'$ contradicts \ref{connectedLeftoverFull_choice:3} since it has fewer components than $G$. Thus, to complete the proof of this claim we need to find an edge of $Z$ in a cycle of $H$. For each $i$, we let $T_{i}$ represent a spanning tree of $C_{i}$. For some $i$, if some $aa'\in E(Z)\cap E(C_{i})$ is not in $T_{i}$, then that edge forms a cycle with edges of $T_{i}$. Consider the situation $E(Z)\cap E(C_{i})\subseteq E(T_{i})$ for all $T_{i}$. If $\delta(G_{0})=1$, then we have the contradiction \[|V|-1\leq |E(Z)|\leq \sum_{i=1}^{t}|E(T_{i})|=\sum_{i=1}^{t}(|C_{i}|-1)\leq|V|-2.\] If $\delta(G_{0})\geq 2$, then given a leaf $a\in V(C_{i})$ there is an edge $aa'\in E(T_{i})\cap E(Z)$ and an edge $ab\notin E(T_{i})$. However, $ab$ forms a cycle with edges of $T_{i}$ that includes $aa'$. Since $H$ is connected, we see by our choice of $G$ that $\delta(H)\geq 2$.
\end{proof}
\end{claim}

We choose a minimally weak set $A_{0}$ of $H$ such that there are critically weak sets $A\subseteq A_{0}$ and $B\subseteq \overline{A_{0}}$ of $H$ with $|\Gamma_{H}(A)|\geq |\Gamma_{H}(B)|$.

\begin{claim}\label{cl:main}If $aa'\in E(Z[A])$ and $bb'\in E(Z[B])$ such that $ab'$ and $ba'$ are non-edges of $G$, then there exists an $X\subset V$ with $\{a,b'\}\subseteq X$ and $\{a',b\}\subseteq \overline{X}$ such that $e_{H}(X,\overline{X})\leq \lambda(H)+2\leq \delta(H)+1$.

\begin{proof}For a contradiction, we suppose $e_{H}(X,\overline{X})\geq \lambda(H)+3$ for every $X\subset V$ with $\{a,b'\}\subseteq X$ and $\{a',b\}\subseteq \overline{X}$. Let $Z'$ be the realization of $\mathcal{R}(Z)$ created by exchanging the edges $aa'$ and $bb'$ in $Z$ with the non-edges $ab'$ and $ba'$ of $G$. Thus, the graph $W=G-E(Z)+E(Z')$ is in $\mathcal{R}(G_{0},G_{0}-E(Z_{0}))$, and by \ref{connectedLeftoverFull_choice:2} we know that $\lambda(W-E(F))=\lambda(H)$. We choose an arbitrary $X\subseteq V$. If either $aa'$ or $b'b$ is not in $E_{H}(X,\overline{X})$, then $e_{W}(X,\overline{X})=e_{H}(X,\overline{X})$. If $\{aa',b'b\}\subseteq E_{H}(X,\overline{X})$, then $e_{W-E(F)}(X,\overline{X})\geq e_{H}(X,\overline{X})-2\geq \lambda(H)+1$. Thus, if $X$ is minimally weak in $W-E(F)$, then it was minimally weak in $H$. Since  $e_{W-E(F)}(A_{0},\overline{A_{0}})=e_{H}(A_{0},\overline{A_{0}})+2$, we know that $A_{0}$ is minimally weak in $H$ and not in $W-E(F)$. Thus, $W-E(F)$ violates \ref{connectedLeftoverFull_choice:3}.
\end{proof}
\end{claim}

\begin{claim}\label{cl:goodFeatures}Let $aa'\in E(Z[A])$ and $bb'\in E(Z[B])$ such that $ab'$ and $ba'$ are non-edges of $G$. If $e_{H}(X,\overline{X})\leq \lambda(H)+2$ for some $X\subseteq V$ with $\{a,b'\}\subseteq X$ and $\{a',b\}\subseteq \overline{X}$, then $\delta(H)$ is odd, $\lambda=\delta(H)-1$, $e_{H}(X,\overline{X})=\delta(H)+1$,
\[e_{H}(X\cap A, \overline{X}\cap A)=e_{H}(X\cap B, \overline{X}\cap B)=\frac{\delta(H)+1}{2}, \text{and}\] 
\[e_{H}(X\cap A, \overline{A})=e_{H}(\overline{X}\cap A, \overline{A})= \frac{\delta(H)-1}{2}.\]
\begin{proof}We have $e_{H}(X,\overline{X})\leq \lambda(H)+2\leq \delta(H)+1$. Applying Lemma~\ref{lem:criticalSet} to both $A$ and $B$ we see that both $e_{H}(X\cap A, \overline{X}\cap A)$ and $e_{H}(X\cap B, \overline{X}\cap B)$ are at least $\bigg\lceil\frac{\delta(H)+1}{2}\bigg\rceil$. Thus, \[\delta(H)+1\geq\lambda(H)+2\geq e_{H}(X,\overline{X})\geq e_{H}(X\cap A, \overline{X}\cap A) + e_{H}(X\cap B, \overline{X}\cap B)\geq 2\bigg\lceil\frac{\delta(H)+1}{2}\bigg\rceil.\] This can only be true if $\delta(H)$ is odd, $\lambda(H)=\delta(H)-1$, $e_{H}(X,\overline{X})=\delta(H)+1$, and 
\[e_{H}(X\cap A, \overline{X}\cap A)=e_{H}(X\cap B, \overline{X}\cap B)=\frac{\delta(H)+1}{2}.\] 
Furthermore, we see by Lemma~\ref{lem:criticalSet} and (\ref{eq:criticalSet3}) that \[\delta(H)-1=e_{H}(A,\overline{A})=e_{H}(X\cap A, \overline{A})+e_{H}(\overline{X}\cap A, \overline{A})\geq \delta(H)-1.\] Thus, \[e_{H}(X\cap A, \overline{A})=e_{H}(\overline{X}\cap A, \overline{A})= \frac{\delta(H)-1}{2}.\qedhere\]
\end{proof}
\end{claim}

\begin{claim}\label{cl:badPath}If there is a path $P$ from a vertex $y\in \Gamma_{H}(A)$ to a $y'\in \Gamma_{H}(B)$ with no internal vertices in $A\cup B$, then either $y$ is adjacent in $F$ to every vertex in $N_{Z}(y')\cap (B-\Gamma_{H}(B))$ or $y'$ is adjacent in $F$ to every vertex in $N_{Z}(y)\cap (A-\Gamma_{H}(A))$.
\begin{proof}Suppose there is an $a\in N_{Z}(y)\cap (A-\Gamma_{H}(A))$ and a vertex $b\in N_{Z}(y')\cap (B-\Gamma_{H}(B))$ such that $a\notin N_{G}(y')$ and $b\notin N_{G}(y)$. Let $X$ be an arbitrary subset of $V$ with $\{a,y'\}\subseteq X$ and $\{b,y\}\subseteq \overline{X}$. By Claim~\ref{cl:goodFeatures} we see that $e_{H}(X\cap A, \overline{X}\cap A)$ and $e_{H}(X\cap B, \overline{X}\cap B)$ are both equal to  $\frac{\delta(H)+1}{2}$. Since $P$ is a path from $y$ to $y'$ that has no internal vertices in $A\cup B$, there must be an edge from $X-B$ to $\overline{X}-A$. Therefore, since $X$ is arbitrary, we see by Claim~\ref{cl:main} the contradiction
\begin{align}
    e_{H}(X,\overline{X})&\geq e_{H}(X\cap A, \overline{X}\cap A)+e_{H}(X\cap B, \overline{X}\cap B
    )+e_{H}(X-B,\overline{X}-A)\notag\\
    &\geq 2\Bigg(\frac{\delta(H)+1}{2}\Bigg)+1\geq \delta(H)+2.\qedhere\notag
    \end{align}
\end{proof}
\end{claim}

We choose an $x\in A$ such that $e_{H}(x,\overline{A})$ is maximized and subject to that we minimize $e_{Z}(x,\overline{A})$.

\begin{claim}\label{cl:notEvenorSmall}$\delta(H)\geq 3$ and odd, $\lambda(H)=\delta(H)-1$, and $e_{H}(x,\overline{A})\leq \frac{\delta(H)-1}{2}$.
\begin{proof}By Lemma~\ref{lem:weakBoundEdge} there is an $a\in A-\Gamma_{H}(A)$ and $b\in B-\Gamma_{H}(B)$. Since $\delta(Z)>\Delta(F)$, we know there exists an $a'\in N_{Z}(a)-N_{F}(b)$ and a $b'\in N_{Z}(b)-N_{F}(a)$. By Claim~\ref{cl:main} there exists an $X\subset V$ with $\{a,b'\}\subseteq X$ and $\{a',b\}\subseteq \overline{X}$ such that $e_{H}(X,\overline{X})\leq \lambda(H)+2\leq \delta(H)+1$. By Claim~\ref{cl:basicCase} and Claim~\ref{cl:goodFeatures},  $\delta(H)$ is at least three and odd, $\lambda=\delta(H)-1$, and \[e_{H}(X\cap A, \overline{A})=e_{H}(\overline{X}\cap A, \overline{A})= \frac{\delta(H)-1}{2}.\] Thus, $e_{H}(x,\overline{A})\leq \frac{\delta(H)-1}{2}$.
\end{proof}
\end{claim}

Since $A$ is weak, we see by Claim~\ref{cl:notEvenorSmall} that $A$ is both minimally and critically weak. We let $A'$ denote all $u\in \Gamma_{H}(A)$ with $e_{H}(u,\overline{A})=\frac{\delta(H)-1}{2}$. Note that $|A'|\leq 2$, and if $A'$ is not empty, then $x\in A'$ by our choice of $x$ and Claim~\ref{cl:notEvenorSmall}. 
 
We say edges $xy$ and $uv$ of a graph $W$ are crossable in $W$ if $E_{G}(\{x,y\},\{u,v\})\neq \emptyset$ and not crossable in $W$, otherwise.
 
\begin{claim}\label{cl:EdgesAreAdjacent}Every $aa'\in E(Z[A])$ and $bb'\in E(Z[B])$ are crossable in $G$.
\begin{proof}By contradiction suppose $aa'$ and $bb'$ are not crossable in $G$. By Claim~\ref{cl:main} there must exist an $X\subset V$ with $\{a,b'\}\subseteq X$ and $\{b,a'\}\subseteq \overline{X}$ and a  $Y\subset V$ with $\{a,b\}\subseteq Y$ and $\{b',a'\}\subseteq \overline{Y}$ such that $e_{H}(X,\overline{X})\leq \delta(H)+1$ and $e_{H}(Y,\overline{Y})\leq \delta(H)+1$. By Claim~\ref{cl:goodFeatures} \[e_{H}(X\cap A, \overline{X}\cap A)=e_{H}(X\cap B, \overline{X}\cap B)=e_{H}(Y\cap A, \overline{Y}\cap A)=e_{H}(Y\cap B, \overline{Y}\cap B)=\frac{\delta(H)+1}{2}, \text{ and}\] \[e_{H}(X\cap A, \overline{A})=e_{H}(X\cap B, \overline{B})=\frac{\delta(H)-1}{2}.\] This implies $e_{H}(X\cap A, \overline{X\cap A})=e_{H}(X\cap A, \overline{X}\cap A)+e_{H}(X\cap A, \overline{A})=\delta(H)$.

Suppose both $Y$ and $\overline{Y}$ intersect $X\cap A$.  Since $\delta(H)$ is odd and $e_{H}(X\cap A, \overline{X\cap A})=\delta(H)$, we see by Lemma~\ref{lem:criticalSet} and (\ref{eq:criticalSet0}) that $e_{H}(Y\cap X\cap A,\overline{Y}\cap X\cap A)\geq \frac{\delta(H)+1}{2}$. However, since $aa'\in E_{H}(Y\cap X\cap A, \overline{Y}\cap \overline{X}\cap A)$, we deduce the contradiction \[e_{H}(Y\cap A, \overline{Y}\cap A)=e_{H}(Y\cap X\cap A,\overline{Y}\cap X\cap A)+E_{H}(Y\cap X\cap A, \overline{Y}\cap \overline{X}\cap A)\geq 1+\frac{\delta(H)+1}{2}.\] An identical argument can be made if $Y$ and $\overline{Y}$ both intersect $X\cap B$. Thus, we assume $X\cap A\subseteq Y$ and $X\cap B\subseteq \overline{Y}$.

Since $H$ is $(\delta(H)-1$)-edge-connected, $e_{H}(X\cap A, \overline{A})=\frac{\delta(H)-1}{2}$, and $\delta(H)\geq 3$, there exists a path $P$ from some $u\in X\cap \Gamma_{H}(A)$ to some $v\in \Gamma_{H}(B)$ such that no internal vertex of $P$ is in $A\cup B$. If $v\in \overline{X}$, then $E(X-B,\overline{X}-A)$ is not empty. However, this leads to the contradiction \[e_{H}(X,\overline{X})\geq e_{H}(X\cap A,\overline{X}\cap A)+e_{H}(X\cap B,\overline{X}\cap B)+e_{H}(X-B,\overline{X}-A)\geq \delta(H)+2.\]  If $v\in X$, then $e_{H}(Y-B,\overline{Y}-A)$ is not empty and we have the contradiction 
\[e_{H}(Y,\overline{Y})\geq e_{H}(Y\cap A,\overline{Y}\cap A)+e_{H}(Y\cap B,\overline{Y}\cap B)+e_{H}(Y-B,\overline{Y}-A)\geq \delta(H)+2.\qedhere\]
\end{proof}
\end{claim}

\begin{claim}\label{cl:FisNotEmpty}$\Delta(F)\geq 1$.
\begin{proof}By contradiction suppose $F$ has no edges, and therefore, $H=G$. We choose an $a\in A-\Gamma_{H}(A)$, $a'\in N_{Z}(a)$, $b\in B-\Gamma_{H}(B)$, and $b'\in N_{z}(b)$. Since $\{b,b'\}\cap N_{G}(a)=\emptyset$ and $\{a,a'\}\cap N_{G}(b)=\emptyset$, we see by Claim~\ref{cl:EdgesAreAdjacent} that $a'b'\in E(G)$. By Claim~\ref{cl:main} there is an $X\subset V$ with $\{a,b'\}\subseteq X$ and $\{b,a'\}\subseteq \overline{X}$ such that $e_{H}(X,\overline{X})\leq \delta(H)+1$. By Claim~\ref{cl:goodFeatures} both $e_{H}(X\cap A,\overline{X}\cap A)$ and $e_{H}(X\cap B,\overline{X}\cap B)$ are equal to $\frac{\delta(H)+1}{2}$. Since $a'b'$ is an edge, we see that $e_{H}(X-B,\overline{X}-A)\geq 1$. However, this gives us the contradiction \[e_{H}(X,\overline{X})\geq e_{H}(X\cap A,\overline{X}\cap A)+e_{H}(X\cap B,\overline{X}\cap B)+e_{H}(X-B,\overline{X}-A)\geq \delta(H)+2.\qedhere\]

\end{proof}
\end{claim}

Since $\Delta(F)\geq 1$, we see that $\delta(Z)\geq 2$.

\begin{claim}\label{cl:speacialS}There exists an $S\subseteq A-\{x\}$ with $|S|\geq 2$ such that $e_{H}(S,A-S)=\frac{\delta(H)+1}{2}$, $e_{H}(S,\overline{A})=\frac{\delta(H)-1}{2}$, and $e_{H}(S,\overline{S})=\delta(H)$.
\begin{proof} 
We choose an $a\in A-\Gamma_{H}(A)$ such that we give priority to those vertices adjacent in $F$ to some vertex in $B-\Gamma_{H}(B)$. Suppose every vertex in $N_{Z}(a)-(A'-\{x\})$ is adjacent in $F$ to every vertex in $B-\Gamma_{H}(B)$. Thus, $\Delta(F)\geq |N_{Z}(a)-(A'-\{x\})|$. However, since $\delta(Z)\geq\Delta(F)+1$ and $|A'|\leq 2$, we see that $|N_{Z}(a)-(A'-\{x\})|=\Delta(F)$, and therefore, $a$ must be adjacent in $Z$ to a vertex in $A'-\{x\}$ and $a$ must not be adjacent in $F$ to a vertex in $B-\Gamma_{H}(B)$. Thus, $A'=\Gamma_{H}(A)$, and since $\Delta(F)\geq 1$, we know $N_{Z}(a)-(A'-\{x\})$ is not empty. Since $|B-\Gamma_{H}(B)|\geq 2$, we see that $\delta(Z)>\Delta(F)\geq 2$ and there must be a vertex in $A-A'$ adjacent to a vertex in $B-\Gamma_{H}(B)$. However, this contradicts our choice of $a$. Thus, there must be some $b\in B-\Gamma_{H}(B)$ not adjacent in $G$ to some $a'\in N_{Z}(a)-(A'-\{x\})$. Let $b'\in N_{Z}(b)-N_{F}(a)$. By Claim~\ref{cl:main} there is an $X\subset V$ with $\{a,b'\}\subseteq X$ and $\{b,a'\}\subseteq \overline{X}$ such that $e_{H}(X,\overline{X})\leq \delta(H)+1$. By Claim~\ref{cl:goodFeatures} $e_{H}(X\cap A, \overline{X}\cap A)=\frac{\delta(H)+1}{2}$ and $e_{H}(X\cap A, \overline{A})=e_{H}(\overline{X}\cap A, \overline{A})=\frac{\delta(H)-1}{2}$. If $x\in \overline{X}$, then since $a\notin \Gamma_{H}(A)$, we see that $|X\cap A|\geq 2$. Suppose $x\in X$. If $a'\notin \Gamma_{H}(A)$, then $|\overline{X}\cap A|\geq 2$, and if $a\in \Gamma_{H}(A)$, then since $a'\notin A'-\{x\}$ and $e_{H}(\overline{X}\cap A, \overline{A})=\frac{\delta(H)-1}{2}$, we see that $|\overline{X}\cap A|\geq 2$. Thus, either $X\cup A$ or $\overline{X}\cup A$ satisfies the conditions of the claim.
\end{proof}
\end{claim}
By Claim~\ref{cl:speacialS} we may choose $S$ to be the smallest such set.

\begin{claim}\label{cl:forbiddenEdge}If there is an edge $aa'\in E(Z[S])$ and a $b\in B-\Gamma_{H}(B)-N_{F}(a')$, then $\{a,a'\}\cap A'\neq \emptyset$.

\begin{proof}Suppose $aa'\in E(Z[S])$ and a vertex $b\in B-\Gamma_{H}(B)-N_{F}(a')$. Let $b'\in N_{Z}(b)-N_{F}(a)$. By Claim~\ref{cl:main} there is an $X\subset V$ with $\{a,b'\}\subseteq X$ and $\{b,a'\}\subseteq \overline{X}$ such that $e_{H}(X,\overline{X})\leq \delta(H)+1$. By Claim~\ref{cl:goodFeatures} $e_{H}(X\cap A, \overline{X}\cap A)=\frac{\delta(H)+1}{2}$ and $e_{H}(X\cap A, \overline{A})=e_{H}(\overline{X}\cap A, \overline{A})=\frac{\delta(H)-1}{2}$. Since $e_{H}(S,\overline{S})=\delta(H)$, we see by Lemma~\ref{lem:criticalSet} and (\ref{eq:criticalSet0}) that $e_{H}(X\cap S, \overline{X}\cap S)\geq \frac{\delta(H)+1}{2}$. This implies $e_{H}(X\cap (A-S),\overline{X}\cap (A-S))=0$ since \[e_{H}(X\cap A, \overline{X}\cap A)=e_{H}(X\cap S, \overline{X}\cap S)+e_{H}(X\cap (A-S),\overline{X}\cap (A-S))\geq \frac{\delta(H)+1}{2}.\]  If $A-S\subset X$, then $\overline{X}\cap A=\overline{X}\cap S$, $|\overline{X}\cap S|<|S|$, $e_{H}(X\cap S, \overline{X}\cap S)=\frac{\delta(H)+1}{2}$,  and $e_{H}(\overline{X}\cap S, \overline{A})=\frac{\delta(H)-1}{2}$. By our choice of $S$, it must be the case $\overline{X}\cap S=\{a'\}$, and therefore, $a'\in A'$. If $A-S\subset \overline{X}$, then $X\cap A=X\cap S$, $|X\cap S|<|S|$, $e_{H}(X\cap S, \overline{X}\cap S)=\frac{\delta(H)+1}{2}$,  and $e_{H}(X\cap S, \overline{A})=\frac{\delta(H)-1}{2}$. By our choice of $S$, it must be the case $X\cap S=\{a\}$, and therefore, $a\in A'$.
\end{proof}
\end{claim}

We choose a $q\in S$ that is adjacent in $H$ to the most vertices in $\overline{A}$. If $A'\cap S\neq \emptyset$, then $\{x,q\}=A'$ by our choice of $x$.

Suppose there is an $a\in S$ such that $N_{Z}(a)\subseteq S$. Since $\delta(Z)\geq 2$, we know that $N_{Z}(a)-\{q\}$ is not empty. Thus, by Claim~\ref{cl:forbiddenEdge} $a$ and every vertex in $N_{Z}(a)-\{q\}$ must be adjacent in $F$ to every $b\in B-\Gamma_{H}(B)$. However, this is a contradiction since $|N_{Z}[a]-\{q\}|\geq \delta(Z)>\Delta(F)$. Therefore, every vertex in $S$ is adjacent in $Z$ to a vertex not in $S$.

If $|S|<\delta(H)$, then \[(\delta(H)-(|S|-1))|S|\leq e_{H}(S,\overline{S})=\delta(H).\] After rearranging and simplifying, we deduce the contradiction $(|S|-1)(\delta(H)-|S|)\leq 0$. Thus, $|S|\geq \delta(H)$, and therefore,  \[\delta(H)=e_{H}(S,\overline{S})\geq e_{Z}(S,\overline{S})\geq |S|\geq \delta(H)\] implies $|S|=\delta(H)$, $E_{H}(S,\overline{S})=E_{Z}(S,\overline{S})$, and every vertex in $S$ is adjacent in $Z$ to exactly one vertex in $\overline{S}$.

Since $|S|=\delta(H)$ and $\delta(Z)\geq 2$, there must be an edge $aa'\in E(Z[S])$ such that $a\in S-N_{F}(b)$ for some $b\in B-\Gamma_{H}(B)$. By Claim~\ref{cl:forbiddenEdge} $a$ or $a'$ must be in $A'$, and therefore, $a$ or $a'$ is $q$. Thus, $A'=\{x,q\}$, and by our choice of $A$ we see that $|\Gamma_{H}(B)|\leq |\Gamma_{H}(A)|=2$. Furthermore, since $q$ is only adjacent to one vertex not in $S$, we see that $\delta(H)=3$. 

Since $|S|=3$ and $q\in A'$, there are two vertices $a$ and $a'$ of $A-\Gamma_{H}(A)$ such that $S=\{q,a,a'\}$. Moreover, Since every vertex in $S$ is adjacent in $Z$ to exactly one vertex in $\overline{S}$ and $\delta(Z)>\Delta(F)\geq 1$, we see by Claim~\ref{cl:forbiddenEdge} that $aq$ and $a'q$ are edges of $Z$. Furthermore, $aa'$ is an edge of $Z$ if and only if $\delta(Z)=3$. To see this, consider that Claim~\ref{cl:forbiddenEdge} says $\Delta(F)=2$ since both $a$ and $a'$ are adjacent in $F$ to every vertex in $\Gamma_{H}(B)$ when $aa'$ is an edge of $Z$. On the other hand, when $\delta(Z)=3$, $aa'$ must be an edge of $Z$.

Since $H$ is $2$ edge-connected and $A'=\{x,q\}$, we know there is a path $P_{q}$ from $q$ to some $q'\in \Gamma_{H}(B)$ and a path $P_{x}$ from $x$ to some $x'\in \Gamma_{H}(B)$ such that $P_{q}$ and $P_{x}$ are edge-disjoint and each have no interval vertices in $A\cup B$. Moreover, if $x'\neq q'$, then the existence of $P_{x}$ and $P_{q}$ lets us reason that $xq'$ and $x'q$ are non-edges of $H$ since $e_{H}(B,\overline{B})<\delta(H)=3$.

Suppose $\Delta(F)=2$. Since $\delta(Z)>\Delta(F)=2$, $\delta(Z)=3$, and therefore, $aa$, $aq$, and $a'q$ are edges of $Z$. Since $|\Gamma_{H}(B)|\leq |\Gamma_{H}(A)|\leq 2$, $e_{H}(B,\overline{B})=\lambda(H)=2<\delta(H)=3$, $\delta(Z)=3$, and $P_{q}$ and $P_{x}$ are edge-disjoint, we may conclude that $q'$ is adjacent in $Z$ to a vertex $b\in B-\Gamma_{H}(B)$. Let $b'\in N_{Z}(b)\cap (B-\Gamma_{H}(B))$. By Claim~\ref{cl:forbiddenEdge} $b$ and $b'$ are adjacent in $F$ to both $a$ and $a'$, and therefore, since $\Delta(F)=2$, $b\notin N_{F}(q)$ and $a\notin N_{F}(q')$. However, this contradicts Claim~\ref{cl:badPath}. 

Thus, we are left with the case $\Delta(F)=1$. We let $u_{0}\in N_{Z}(a)\cap (A-S)$ and $u_{1}\in N_{Z}(a')\cap (A-S)$.

Suppose $q'$ is adjacent in $Z$ to a vertex $b\in B-\Gamma_{H}(B)$. Let $b'\in N_{Z}(b)\cap (B-\Gamma_{H}(B))$. Since both $a$ and $a'$ cannot be adjacent in $F$ to $q'$, Claim~\ref{cl:badPath} implies $q$ is adjacent in $F$ to $b$.  By Claim~\ref{cl:EdgesAreAdjacent} both $b'$ and $q'$ are adjacent in $G$ to at least one end of each of the edges $au_{0}$ and $a'u_{1}$. Thus, one of $b'$ and $q'$ are adjacent in $G$ to two vertices in $\{a,a',u_{0},u_{1}\}$. Recall that $\{u_{0},u_{1}\}\subseteq A-S$, $\Gamma_{H}(A)=\{x,q\}$, and $x$ could be one or both of $u_{0}$ and $u_{1}$. However, since $xq'$ is a non-edge of $H$, any edge from  $b'$ and $q'$ to $\{a,a',u_{0},u_{1}\}$ must be an edge of $F$. This contradicts $\Delta(F)=1$. Thus, $q'$ is not adjacent in $Z$ to a vertex in $B-\Gamma(B)$. Since $\delta(Z)\geq 2$, there must be a cycle $C$ in $Z$ such that $V(C)\subseteq B-\{q'\}$. 

Since $x'q$ is a non-edge of $H$ and $V(C)-\{x'\}\subseteq \Gamma_{H}(B)$, we may deduce that $a$, $a'$, and $q$ are not adjacent in $H$ to vertices in $C$. Let $v\in \{a,a',q,u_{0},u_{1}\}$. If there is a $v'\in N_{Z}(v)\cap A$ that is not adjacent in $H$ to a vertex of $C$, then we have a contradiction since Claim~\ref{cl:EdgesAreAdjacent} implies that $v$ must be adjacent in $F$ to at least two vertices on $V(C)$. Thus, every vertex in $N_{Z}(v)$ is adjacent in $H$ to a vertex in $V(C)$. Suppose $|C|\geq 4$, and let $v_{0}v_{1}$, $v_{1}v_{2}$, and $v_{2}v_{3}$ be consecutive edges along $C$ such that $q$ is adjacent in $F$ to $v_{0}$. Thus, $q$ is not adjacent in $F$ to $\{v_{1},v_{2},v_{3}\}$. By Claim~\ref{cl:EdgesAreAdjacent} either $a$ or $a'$ must be adjacent in $F$ to at least two vertices in $\{v_{1},v_{2},v_{3}\}$. This contradiction implies $|C|=3$. Since $\{a,a',q\}$ are all adjacent in $F$ to distinct vertices in $C$ and $u_{0}$ and $u_{1}$ are adjacent in $H$ to vertices in $C$ we may conclude that $u_{0}=u_{1}=x$ and $x$ must be adjacent in $H$ to $x'$. This implies $x'$ is adjacent along $C$ to two vertices in $B-\Gamma_{H}(B)$. We have a contradiction since Claim~\ref{cl:badPath} says either $x'$ is adjacent in $F$ to both $a$ and $a'$ or $x$ is adjacent in $F$ to two vertices on $C$. This completes the proof of Theorem~\ref{thm:connectedLeftoverFull}.
\end{proof}

% \bibliographystyle{amsplain}
% \bibliography{Bibliography/main}

\begin{thebibliography}{10}

\bibitem{Berge1958}
Claude Berge, \emph{Th{\'e}orie des graphes et ses applications}, 1 ed.,
  Collection universitaire de math{\'e}matiques, Dunod, Paris, 1958.

\bibitem{Bollobas1979}
B{\'{e}}la Bollob{\'{a}}s, \emph{On graphs with equal edge connectivity and
  minimum degree}, Discrete Mathematics \textbf{28} (1979), no.~3, 321--323.

\bibitem{Bollobas1985}
B.~Bollobás, Akira Saito, and N.~C. Wormald, \emph{Regular factors of regular
  graphs}, Journal of Graph Theory \textbf{9} (1985), no.~1, 97--103.

\bibitem{Brualdi1978}
R.~A. Brualdi, \emph{Probl{\'{e}}mes}, Probl{\'{e}}mes combinatoires et
  th{\'{e}}orie des graphes, Colloq. Internat., vol. 260, Paris, 1978,
  pp.~437--443.

\bibitem{Busch2012}
Arthur~H. Busch, Michael~J. Ferrara, Stephen~G. Hartke, Michael~S. Jacobson,
  Hemanshu Kaul, and Douglas~B. West, \emph{Packing of graphic n-tuples},
  Journal of Graph Theory \textbf{70} (2012), no.~1, 29--39.

\bibitem{Chetwynd1985}
A.~G. Chetwynd and A.~J.~W. Hilton, \emph{Regular graphs of high degree are
  1-factorizable}, Proceedings of the London Mathematical Society
  \textbf{s3-50} (1985), no.~2, 193--206.

\bibitem{Csaba2016}
B{\'{e}}la Csaba, Daniela Kühn, Allan Lo, Deryk Osthus, and Andrew Treglown,
  \emph{Proof of the 1-factorization and {H}amilton decomposition conjectures},
  vol. 244, American Mathematical Society ({AMS}), nov 2016.

\bibitem{Dankelmann2000}
Peter Dankelmann and Lutz Volkmann, \emph{Degree sequence conditions for
  maximally edge-connected graphs depending on the clique number}, Discrete
  Mathematics \textbf{211} (2000), no.~1, 217 -- 223.

\bibitem{Diestel2016}
Reinhart Diestel, \emph{Graph theory}, fifth ed., Graduate Texts in
  Mathematics, vol. 173, Springer-Verlag, Heidelberg, August 2016.

\bibitem{Edmonds1964}
J.~Edmonds, \emph{Existence of k-edge connected ordinary graphs with prescribed
  degrees}, Journal of Research of the National Bureau of Standards Section B
  Mathematics and Mathematical Physics \textbf{68B} (1964), no.~2, 73 -- 74.

\bibitem{Goldsmith1978}
Donald~L. Goldsmith and Arthur~T. White, \emph{On graphs with equal
  edge-connectivity and minimum degree}, Discrete Mathematics \textbf{23}
  (1978), no.~1, 31 -- 36.

\bibitem{Gu2013}
Xiaofeng Gu and Hong-Jian Lai, \emph{Realizing degree sequences with
  k-edge-connected uniform hypergraphs}, Discrete Mathematics \textbf{313}
  (2013), no.~12, 1394--1400.

\bibitem{Hellwig2003}
Angelika Hellwig and Lutz Volkmann, \emph{Maximally edge-connected digraphs},
  Australasian Journal of Combinatorics \textbf{27} (2003), 23--32.

\bibitem{Hellwig2008}
Angelika Hellwig and Lutz Volkmann, \emph{Maximally edge-connected and
  vertex-connected graphs and digraphs: A survey}, Discrete Mathematics
  \textbf{308} (2008), no.~15, 3265 -- 3296.

\bibitem{Katerinis1993}
P.~Katerinis, \emph{Regular factors in regular graphs}, Discrete Mathematics
  \textbf{113} (1993), no.~1, 269--274.

\bibitem{Kundu1974}
Sukhamay Kundu, \emph{Generalizations of the k-factor theorem}, Discrete
  Mathematics \textbf{9} (1974), no.~2, 173--179.

\bibitem{Mattiolo2022}
Davide Mattiolo and Eckhard Steffen, \emph{Highly edge-connected regular graphs
  without large factorizable subgraphs}, Journal of Graph Theory \textbf{99}
  (2022), no.~1, 107--116.

\bibitem{Plesnik1972}
J{\'a}n Plesn{\'\i}k, \emph{Connectivity of regular graphs and the existence of
  1-factors}, Matematick{\`y} {\v{c}}asopis \textbf{22} (1972), no.~4,
  310--318.

\bibitem{Seacrest2021}
Tyler Seacrest, \emph{Multi-switch: a tool for finding potential edge-disjoint
  $1$-factors}, Electronic Journal of Graph Theory and Applications \textbf{9}
  (2021), no.~1, 87--94.

\bibitem{Shiu2008}
Wai~Chee Shiu and Gui~Zhen Liu, \emph{k-factors in regular graphs}, Acta
  Mathematica Sinica, English Series \textbf{24} (2008), no.~7, 1213--1220.

\bibitem{Thomassen2020}
Carsten Thomassen, \emph{Factorizing regular graphs}, Journal of Combinatorial
  Theory, Series B \textbf{141} (2020), 343--351.

\bibitem{Tian2014}
Yingzhi Tian, Jixiang Meng, Hongjian Lai, and Zhao Zhang, \emph{On the
  existence of super edge-connected graphs with prescribed degrees}, Discrete
  Mathematics \textbf{328} (2014), 36--41.

\bibitem{Wang1974}
D.~L. Wang and D.~J. Kleitman, \emph{A note on n-edge-connectivity}, SIAM
  Journal on Applied Mathematics \textbf{26} (1974), no.~2, 313--314.

\end{thebibliography}

\providecommand{\bysame}{\leavevmode\hbox to3em{\hrulefill}\thinspace}
\providecommand{\MR}{\relax\ifhmode\unskip\space\fi MR }
\providecommand{\MRhref}[2]{%
  \href{http://www.ams.org/mathscinet-getitem?mr=#1}{#2}
}
\providecommand{\href}[2]{#2}

\end{document}